\input amstex
\documentstyle{amsppt}

\topmatter
\title Complete ccc Boolean algebras, the order
sequential topology, and a problem of von Neumann \endtitle

\author B.Balcar, T.Jech and T.Paz\' ak \endauthor

\affil
Mathematical Institute, Academy of Sciences of the Czech Republic
\endaffil

\thanks Supported in part by the GA\v CR Grant number 201/02/0857 (Balcar and 
Jech)
and the GA\v CR Grant number 201/03/0933 (Paz\'ak).
\endthanks

\dedicatory
On the occasion of John von Neumann's 100th birthday
\enddedicatory

\date December 28, 2003 \enddate

\address
\v Zitn\'a 25, 115~67~Praha 1, Czech Republic
\endaddress

\email balcar\@cts.cuni.cz, jech\@math.cas.cz,
tpazak\@yahoo.com \endemail

\keywords Boolean algebra, measure, submeasure, weak
distributivity, independent reals, sequential topology 
\endkeywords

\rightheadtext{Complete ccc Boolean algebras}
% the sequential topology, and a problem of von Neumann}
\leftheadtext{B.Balcar, T.Jech and T.Paz\' ak}

\abstract Let $B$ be a complete $ccc$ Boolean algebra and let
$\tau_s$ be the topology on $B$ induced by the algebraic
convergence of sequences in $B$.
\roster
\item "1." Either there 
exists a Maharam submeasure on $B$ or every nonempty open set 
in $(B,\tau_s)$ is topologically dense.
\item "2." It is consistent that every weakly distributive
complete $ccc$ Boolean algebra carries a strictly positive
Maharam submeasure.
\item "3." The topological space $(B,\tau_s)$ is sequentially
compact if and only if the generic extension by $B$ does not
add independent reals.
\endroster
We also give examples of $ccc$ forcings adding a real
but not independent reals.
\endabstract
\endtopmatter

%.............................endtopmatter................

\document

\subheading{1. Introduction} We investigate combinatorial
properties of
complete $ccc$ Bool\-ean algebras. The focus is on properties
related to the existence of a Maharam submeasure and on 
forcing properties. In particular, we address the question
when the forcing adds independent reals. 
The work is continuation of \cite{BGJ} and \cite{BFH} and
is related to 
the problems of 
von~Neumann and Maharam. 

The problem of von~Neumann
from The Scottish Book (\cite{Sc}, Problem 163) asks whether every
weakly distributive complete $ccc$ Boolean algebra
carries a countably additive measure. 

Von Neumann's problem can be divided into two distinctly
different questions.
Weak distributivity
is a consequence of a property possibly weaker
then measurability, namely the existence of a continuous 
strictly positive submeasure (a Maharam submeasure); the
Control Measure Problem of \cite{M} asks whether every complete
Boolean algebra that carries a continuous submeasure must also
carry a measure. For exact formulation of this question see \cite{F}.
It should be noted that the Control Measure Problem is eqiuvalent
to a $\Pi^1_2$ statement and is therefore absolute for inner models
and generic extensions.

The second question is this modified von Neumann's problem:
does every weakly distributive $ccc$ complete Boolean algebra
carry a strictly positive Maharam submeasure? This statement
is not provable in ZFC, as the algebra associated with a Suslin
tree is a counterexample; moreover, a forcing notion constructed by Jensen
under $V=L$ in \cite{J} yields a counterexample that is countably
generated. In Section 3 we show that it is consistent that the
modified von Neumann's problem holds.

D.~Maharam \cite{M} characterized algebras that carry a continuous 
submeasure as those on which the sequential topology
$\tau_s$ is metrizable. In \cite{BGJ} this is improved
to the condition that $B$ is $ccc$ and $(B,\tau_s)$ is
a Hausdorff space. In Section~2 we prove the following

\proclaim{Decomposition Theorem} Let $B$ be a 
complete $ccc$ Boolean algebra. Then 
there are disjoint elements $d,m \in B$ such that
$d \vee m = \bold 1$ and
\roster
\item "(i)" In the space $(B \upharpoonright d,\tau_s)$
the closure of every nonempty open set is the whole space.
\item "(ii)" The Boolean algebra $B \upharpoonright m$
carries a strictly positive Maharam submeasure.
\endroster
\endproclaim

A sequence $\langle a_n : n \in \omega \rangle$ in a 
topological space {\it converges} to a point $a$ ($\lim_n a_n = a$)
if for every open neighborhood $U$ of $a$, all but 
finitely many $a_n$'s belong to $U$. If the space is
Hausdorff then every sequence converges to at most one 
point (the {\it limit} of $a_n$). A topological space $X$ is 
{\it Fr\' echet} if for every set $A \subset X$, every point
in the closure of $A$ is the limit of some sequence
in $A$. A space is {\it sequentially compact} if every
sequence has a convergent subsequence.

Let $B$ be a complete Boolean algebra. $B$ is 
{\it weakly distributive} (more exactly $(\omega,
\omega)$-weakly distributive) if for every sequence
$\langle P_n : n \in \omega \rangle$ of countable
maximal antichains, there exists a maximal antichain $Q$
with the property that each $q \in Q$ meets only
finitely many elements of each $P_n$. Equivalently,
$B$ is weakly distributive if and only if in any
generic extension $V[G]$ by $B$, every function
$f:\omega \rightarrow \omega$ is bounded by some
$g: \omega \rightarrow \omega$ in the ground model $V$
(i.e. $f(n) \leq g(n)$ for all $n \in \omega$).

$B$ {\it adds an independent real} if there exists 
some $X \subset \omega$ in 
$V[G]$, the generic extension by $B$, such that neither $X$ nor
its complement has an infinite subset $Y$ such that
$Y \in V$. Note that both Cohen and random forcing add independent
reals; the latter is weakly distributive while the former
is not. Neither Sacks forcing \cite{S} (weakly distributive)
nor Miller forcing \cite{Mi} (not weakly distributive) add independent
reals. In Section~5 we present $ccc$ variants of these forcings.

A sequence $\langle a_n : n \in \omega \rangle$ in $B$ 
{\it converges} (algebraically) to $\bold 0$ if there exists
a decreasing sequence $b_0 \geq b_1 \geq \dots \geq b_n \geq \dots$
with $\bigwedge_{n \in \omega} b_n = \bold 0$ such that 
$a_n \leq b_n$ for all $n \in \omega$. A sequence 
$\langle a_n : n \in \omega \rangle$ converges to $a \in B$
($\lim a_n = a$) if the sequence $\langle a_n \vartriangle a 
: n \in \omega \rangle$ of symmetric differences converges
to $\bold 0$. The {\it order sequential topology} $\tau_s$ on
$B$ is defined as follows: the closure of a set $A \subset B$
is the  smallest set $\overline A \supset A$ with the property
that the  limit of every convergent sequence in $\overline A$
is in
$\overline A$. The space $(B,\tau_s)$ is $T_1$ and every 
topologically convergent sequence has a unique limit. Moreover,
a sequence $\langle a_n : n \in \omega \rangle$ in $B$ converges to
$a$ topologically if and only if every subsequence of
$\langle a_n : n \in \omega \rangle$ has a subsequence
that converges to $a$ algebraically.

In \cite{BGJ} it is proved that the space $(B,\tau_s)$
is Fr\' echet if and only if $B$ is weakly distributive
and satisfies the $\goth b$-chain condition, where $\goth b$
is the {\it bounding number}: the least cardinality
of a family $\Cal F$ of functions from $\omega$ to 
$\omega$ such that $\Cal F$ is unbounded; i.e. for
every $g: \omega \rightarrow \omega$ 
there is some $f \in \Cal F$
such that $g(n) \leq f(n)$ for infinitely many $n$'s.

A {\it submeasure} on $B$ is a nonnegative real valued 
function $\mu$ such that
\roster
\item "(i)" $\mu(\bold 0) = 0$
\item "(ii)" $\mu(a) \leq \mu(b)$ whenever $a \leq b$
\item "(iii)" $\mu(a \vee b) \leq \mu(a) + \mu(b)$
\endroster

A {\it strictly positive submeasure} has 
$\mu(a) = 0$ only if $a=\bold 0$.
A {\it Maharam submeasure} is continuous; 
i.e. $\lim \mu (a_n) = 0$ for every decreasing sequence 
$\langle a_n : n \in \omega \rangle$ such that 
$\bigwedge_n a_n = \bold 0$.

If $B$ carries a strictly positive Maharam submeasure then
$B$ is $ccc$ and weakly distributive, and if atomless then
it adds an independent real.

By \cite{M} and \cite{BGJ}, $B$ carries strictly positive
Maharam submeasure if and only if $(B,\tau_s)$ is metrizable,
if and only if $(B,\tau_s)$ is regular, if and only if $B$
is $ccc$ and $(B,\tau_s)$ is Hausdorff.

%..............2. decomposition theorem......................................

\subheading{2. Decomposition theorem} We shall now use the results and 
techniques from \cite{BGJ} to prove the Decomposition theorem.

\proclaim{Theorem} Let $B$ be a 
complete $ccc$ Boolean algebra. Then 
there are disjoint elements $d,m \in B$ such that
$d \vee m = \bold 1$ and
\roster
\item "(i)" In the space $(B \upharpoonright d,\tau_s)$
the closure of every nonempty open set is the whole space.
\item "(ii)" The Boolean algebra $B \upharpoonright m$
carries a strictly positive Maharam submeasure.
\endroster
\endproclaim

The elements $d$ and $m$ are uniquely determined, and either can be 
$\bold 0$. If $m \neq \bold 0$ then $B$ carries a nontrivial
continuous submeasure while if $m = \bold 0$ then every
continuous real valued function on $B$ is constant.

\demo{Proof} First we prove the theorem in the case when $B$
is weakly distributive.
Let $B$ be a weakly distributive complete $ccc$ Boolean algebra.
By \cite{BGJ} the space $(B,\tau_s)$ is Fr\' echet. Let
$$\Cal N = \{U: U \ \text{is an open neighborhood of}\ \bold 0 
\ \text{and is downward closed} \}$$
where downward closed means that $a < b \in U$ implies $a \in U$.
It is proved in \cite{BGJ} (Lemma 3.6) that if $(B,\tau_s)$
is Fr\' echet then $\Cal N$ is a neighborhood base of $\bold 0$,
that $U \vee V = U \vartriangle V$ for $U,V \in \Cal N$, and
that the closure $cl(A)$ of each downward closed set $A$ is
$\bigcap \{A \vee U : U \in \Cal N \}$ and is also downward closed.

Now let $$D = \bigcap\{cl(U) : U \in \Cal N \} \ \text{and} \ 
d = \bigvee D.$$ $D$ is both downward closed and topologically
closed, and it follows from the remarks above that
$$D = \bigcap \{U \vee V : U,V \in \Cal N \} = 
\bigcap \{U \vee U : U \in \Cal N \}.$$
If $a \not \in D$ then for some $U \in \Cal N$, $a \not \in U \vee U$
and hence $U$ and $a \vartriangle U$ are disjoint; in other
words, $a$ is Hausdorff separated from $\bold 0$.
It follows that if we let $m = -d$, then the space 
$(B\upharpoonright m,\tau_s)$ is a Hausdorff space. By \cite{BGJ},
$B\upharpoonright m$ carries a strictly positive Maharam submeasure.
It remains to show that in $(B\upharpoonright d,\tau_s)$, every
nonempty open set is dense.
\enddemo %just virtual end

\proclaim{Lemma} $D$ is closed under $\vee$.\endproclaim

\demo{Proof} Let $U \in \Cal N$. By \cite{BGJ} (Lemma 4.11),
there exists a $U_1 \in \Cal N$ such that $U_1 \vee U_1 \vee U_1
\subset U \vee U$. Similarly, there exists a $U_2 \in \Cal N$ 
such that $U_2 \vee U_2 \vee U_2 \subset U_1 \vee U_1$, and
letting $V = U_1 \cap U_2$ we get $V \vee V \vee  V \vee V 
\subset U \vee U$. Hence 
$D = \bigcap \{ V \vee V \vee  V \vee V : V \in \Cal N \}$.
It follows that $D \vee D = D$.
\hfill $\square$
\enddemo

Now let $\{a_n : n \in \omega \}$ be a maximal antichain in $D$.
The sequence $\{ \bigvee_{k=0}^n a_k :n \in \omega \}$ is in $D$
and converges to $d$. Since $D$ is closed, we have $d \in D$, and
so $B \upharpoonright d = D$.

For every $U \in \Cal N$, $cl(U)=\{ U \vee V : V \in \Cal N \}$
and hence $cl(U) \supseteq D$. Now let $G$ be an arbitrary
topologically open set in $B \upharpoonright d$. There exist
$a \in D$ and $U \in \Cal N$ such that $G \supseteq (a \vartriangle U)
\cap D$. Since $cl(a \vartriangle U) \supseteq a \vartriangle D = D$,
we have $cl(G) \supseteq D$ and the theorem follows for the weakly
distributive case.

\smallskip

In the general case, there exists an element $d_1 \in B$ such that
$B \upharpoonright -d_1$ is weakly distributive, and such that
$B \upharpoonright d_1$ is nowhere weakly distributive. There
exists an infinite matrix $\{a_{kl} \}$ such that each row is
a partition of $d_1$ and for every nonzero $x \leq d_1$ there is
some $k \in \omega$ such that $x \wedge a_{kl} \neq \bold 0$
for infinitely many $l$.

Let $d_2$ and $m$, with $d_2 \vee m = -d_1$, be the decomposition
of the weakly distributive algebra $B \upharpoonright -d_1$, so
that $B \upharpoonright m$ carries a strictly positive Maharam 
submeasure and $(B \upharpoonright d_2,\tau_s)$ has the property
that every nonempty open set is dense in the space. Let
$d= d_1 \vee d_2$, and let us prove that in 
 $(B \upharpoonright d,\tau_s)$, every nonempty open set is dense.

Let $U$ be an open neighborhood of $\bold 0$ in 
$B \upharpoonright d$. The space $(B \upharpoonright d_2,\tau_s)$
is a closed subspace of $(B \upharpoonright d,\tau_s)$, and
$V = U \cap B \upharpoonright d_2$ is an open neighborhood of $\bold 0$ in 
$B \upharpoonright d_2$. 

Let $c \leq d$ be arbitrary; we shall prove that $c$ is in the
closure of $U$. Let $c_1 = c \wedge d_1$ and $c_2 = c \wedge d_2$.
Since $c_2$ is in the closure of $V$ and $B \upharpoonright d_2$
is Fr\' echet, there exists a sequence 
$\langle z_n : n \in \omega \rangle$ in $V$ that converges to $c_2$.
We shall prove that $c_1 \vee z_n \in cl(U)$ for each $n \in \omega$,
and then it follows that $c = \lim_n (c_1 \vee z_n)$ is in $cl(U)$.

Thus let $n \in \omega$ be fixed. For every $k$ and every $l$ let
$y_{kl} = c_1 \wedge \bigvee_{i\geq l} a_{ki}$. Since the sequence
$\langle y_{0l} : l \in \omega \rangle$ converges to $\bold 0$,
we have  $\lim_l (y_{0l}\vee z_n) = z_n$, and since $z_n \in U$,
there exists some $l_0$ such that $y_{0l_0}\vee z_n \in U$.
Let $x_0 = y_{0l_0}$.

Next we consider the sequence $\langle y_{1l}\vee x_0 \vee z_n :
l \in \omega \rangle$. This sequence converges to $x_0 \vee z_n
\in U$ and so there exists some $l_1$ such that $x_1 \vee z_n \in U$
where $x_1 = y_{1l_1}\vee x_0$. We proceed by induction and obtain
a sequence $\langle l_k : k \in \omega \rangle$ and an increasing sequence 
$\langle x_k : k \in \omega \rangle$ with
$x_k \vee z_n \in U$ for each $k$. The sequence $\langle x_k : 
k \in \omega \rangle$ converges to $c_1$ because otherwise, if
$b \neq \bold 0$ is the complement of $\bigvee_k x_k$ in $c$,
then
$b \leq \bigwedge_k \bigvee_{i < l_k} a_{ki}$ and so $b$ meets
only finitely many elements in each row of the matrix. Hence
$c_1 \vee z_n = \lim_k (x_k \vee z_n) \in cl(U)$.
\hfill $\square$

\smallskip

Let us note that in the case $d=\bold 1$ we have $U \vee U = B$
for every open neighborhood
$U \in \Cal N$.

%.........................3. weak distributivity......................

\subheading{3. Weak distributivity} If $B$ is a complete
$ccc$ Boolean algebra then $B$ is weakly distributive if
and only if the space $(B,\tau_s)$ is Fr\' echet. In
this Section we present yet another necessary and sufficient
condition for weak distributivity.

\definition{Definition} $I_s$ is the collection of all sets 
$A \subset B$ such that $A$ is either finite, or
is the range of a sequence in $B$
that converges to $\bold 0$.
\enddefinition

$I_s$ is an ideal of sets, $I_s \subset \Cal P (B)$. An ideal
$I$ of sets is a {\it p-ideal} if for any sequence $\langle A_n 
: n \in \omega \rangle$ of sets in $I$ there exists a set $A \in I$
such that $A_n - A$ is finite, for every $n \in \omega$. We
shall prove the following equivalence (one direction was obtained
independently by S.~Quickert \cite{Q}):

\proclaim{Theorem} A complete $ccc$ Boolean algebra $B$ is weakly 
distributive if and only if the ideal $I_s$ is a p-ideal.
\endproclaim

First we give a different description of $I_s$:

\proclaim{Lemma} Let $B$ be a complete $ccc$ Boolean algebra. 
Then $A \in I_s$ if and only if there exists a maximal 
antichain $W$ such that every $w \in W$ is incompatible with
all but finitely many $a \in A$.
\endproclaim

\demo{Proof} First let $A \in I_s$, $A =\{a_n : n \in \omega\}$
where $\lim a_n = \bold 0$. Let $b_n \geq a_n$ be such that
$\langle b_n : n \in \omega \rangle$ is decreasing and 
$\bigwedge_n b_n = \bold 0$. We may assume that $b_0 = \bold 1$,
and set $w_n = b_n - b_{n+1}$, for each $n \in \omega$. The set
$W = \{w_n : n \in \omega \}$ is a maximal antichain and each 
$w_n$ is incompatible with all $a_k$, $k \geq n+1$.

Conversely let $A$ satisfy the condition of the lemma, with
$W = \{w_n : n \in \omega \}$ an antichain that witnesses it.
If $A$ is infinite then it is necessarily countable, say 
$A = \{a_n : n \in \omega \}$. We claim that 
$\langle a_n : n \in \omega \rangle$ converges to $\bold 0$.
For each $n$, let $b_n = \bigvee_{k=n}^\infty w_k$; the sequence
$\langle b_n : n \in \omega \rangle$ is decreasing and 
$\bigwedge_n b_n = \bold 0$. From the condition on $W$ it
follows that there is an increasing sequence $k_0 < k_1 < \dots
< k_n \dots$ such that for every $n \in \omega$, $b_n \geq a_k$
for all $k \geq k_n$. This implies that $\langle a_n : n \in 
\omega \rangle$ converges to $\bold 0$.
\hfill $\square$
\enddemo

\demo{Proof of the theorem} First let us assume that $B$ is
weakly distributive, and let $\langle A_n : n \in \omega \rangle$
be a sequence of sets in $I_s$. For each $n \in \omega$, let
$W_n$ be a witness to $A_n \in I_s$. By weak distributivity there is 
a maximal antichain $W = \{w_n : n \in \omega \}$ such that
for each $n$ and each $k$, $w_n$ meets only finitely many
elements of $W_k$. It follows that for each $n$ and each $k$,
$w_n$ meets only finitely many elements of $A_k$. Let 

$$A = \{a: \exists n \ a \in A_n \ \text{and} \ a \
\text{is incompatible with every}\ w_i,\ i \leq n \}.$$
For any given $n$, if $a \in A_n - A$ then $a$ is compatible with
some $w_i$, $i\leq n$; since for each $i$ there are only 
finitely many such $a$, the set $A_n - A$ is finite. To complete
the proof we show that $W$ witnesses that $A \in I_s$:
For each $i$, if $w_i$ meets $a \in A$ then $a \not \in A_n$ for
all $n \geq i$ and hence $a \in A_k$ for some $k<i$.
Therefore $w_i$ meets only finitely many $a \in A$.

Conversely, assume that $I_s$ is a p-ideal, and let
$\langle W_n : n \in \omega \rangle$ be a sequence of maximal 
antichains.
Since every maximal antichain is itself in $I_s$, there exists
a set $A \in I_s$ such that $W_n - A$ is finite for every $n \in 
\omega$.
Let $W$ be a witness to $A \in I_s$. Each $w \in W$ meets only 
finitely many $a \in A$, and since $W_n - A$ is finite, $w$ meets 
only finitely many elements of $W_n$, for each $n$.
This proves that $B$ is weakly distributive.
\hfill $\square$
\enddemo

Let us recall Todor\v cevi\' c's dichotomy for 
p-ideals of countable sets  usually denoted by (*):

\proclaim{(*) Dichotomy for p-ideals}\cite{T}
Let $S$ be an infinite set. Then for every p-ideal $\Cal I \subset
[S]^{\leq \omega}$ either 
\roster
\item"(i)" $\exists Y \subset S$ uncountable such that $[Y]^{\leq
\omega} \subset  \Cal I$, or
\item"(ii)" $\exists \{S_n : n \in \omega \}$ such that $\bigcup_n S_n
= S \quad \text{and} \quad \forall n \in \omega \quad 
\forall I \in \Cal I \quad |S_n \cap I | < \omega$.
\endroster
\endproclaim

As it is shown in \cite{T} the dichotomy is a consequence of the 
Proper Forcing Axiom as well as it is consistent with GCH relative 
to the consistency of the existence of a supercompact cardinal
(\cite{AT}).

\proclaim{Lemma} If $B$ is $ccc$ then there
exists no uncountable $X \subset B$ such that 
$[X]^\omega \subset I_s$.
\endproclaim

\demo{Proof} Let $X \subset B$ be uncountable.
Since $B$ is $ccc$ there is some $b \in B^+$ such that
any $\bold 0 < a \leq b$ meets uncountably many elements from
$X$. Let $X_0 \in [X]^\omega$ be such that $\bigvee X_0 \geq b$
and by induction let $X_n \in [X - \bigcup_{i<n}X_n]^\omega$ 
such that $X_n \geq b$. There is no antichain $W$ witnessing 
$\bigcup X_n \in I_s$ since any $w \in W$ compatible with
$b$ is also compatible with infinitely many elements of 
$\bigcup X$. \hfill $\square$
\enddemo

\proclaim{Corollary} Let $B$ be weakly distributive, $ccc$ complete
Boolean algebra. (*) implies that every singleton is a $G_\delta$
set in $(B,\tau_s)$.
\endproclaim

\demo{Proof} It is enough to show that $\{\bold 0\}$ is a $G_\delta$
set in $(B,\tau_s)$. Assuming (*) for $I_s$, it follows that 
$B = \bigcup_{n=0}^\infty S_n$ with each $S_n$ meeting only 
finitely many elements of each $A \in I_s$. It follows that $\bold 0$ is
not in the closure of $S_n - \{\bold 0 \}$ for any $n$.
%Note that $(B,\tau_s)$ is Fr\' echet since $B$ is $ccc$ weakly
%distributive Boolean algebra. 
Let $U_n = B - \text{cl}(S_n -\{\bold 0 \})$. Each $U_n$ is
an open neighborhood of $\bold 0$ and 
$\bigcap_{n \in \omega} U_n = \{\bold 0\}$, hence $\{\bold 0 \}$ 
is $G_\delta$ in $(B,\tau_s)$.
\hfill $\square$
\enddemo

\proclaim{Theorem} Assuming (*), every weakly distributive
$ccc$ complete Boolean algebra carries a strictly positive Maharam
submeasure.
\endproclaim

\demo{Proof} Let $m,d \in B$ be given by the decomposition theorem.
If $m = \bold 1$ the space $(B,\tau_s)$ is completely metrizable.
Suppose now that $d > \bold 0$. By the previous corollary there is
a family $\{U_n : n \in \omega \}$ of open neighborhoods of $\bold 0$
such that $\bigcap_{n \in \omega} U_n = \{\bold 0\}$. We may assume that
$U_{n+1}\subset U_n$, and since $B$ is
weakly distributive, the space is Fr\'echet and we may assume that
each $U_n$ is downward closed. By the decomposition theorem,
$d$ is in the closure of every nonempty open set, and since the
space is Fr\'echet, there exists for each $n$ a sequence $\{a^n_k\}_k$
in $U_n$ that converges to $d$. By weak distributivity there exists
a function $k(n)$ such that the sequence $\{b_n\}_n=\{a^n_{k(n)}\}_n$ 
converges to $d$. Since $d>\bold 0$ there exists a $c>\bold 0$ such that
$b_n\ge c$ for eventually all $n$, say all $n\ge n_0.$. 
Since each $U_n$ is downward closed
and $b_n\in U_n$, it follows that $c\in U_n$ for all $n\ge n_0$, a contradiction.
\hfill $\square$
\enddemo

%......................4. independent reals.........................

\subheading{4. Independent reals}We shall now give  
a necessary and sufficient condition for a complete $ccc$ 
Boolean algebra to add independent reals.

\proclaim{Theorem} Let $B$ be a complete $ccc$ Boolean algebra.
$B$ does not add independent reals if and only if $(B,\tau_s)$
is sequentially compact.
\endproclaim

\demo{Proof} We identify infinite sequences $\langle a_n : 
n \in \omega \rangle$ in $B$ with Boolean names for subsets
of $\omega$, namely $a_n = |\! | n \in \dot A |\! |$

First let $\langle a_n : n \in \omega \rangle$ be a name for
an independent real (without loss of generality we assume that
it is independent with Boolean value $\bold 1$). We shall prove
that $\langle a_n : n \in \omega \rangle$ has no topologically
convergent subsequence. Toward a contradiction, assume that
it does. Then it has an algebraically convergent subsequence
and since an independent real intersected with an infinite
ground model set is independent on that set, we may as well
assume that $\langle a_n : n \in \omega \rangle$ itself
is convergent. Let $a$ be the limit of $\langle a_n : 
n \in \omega \rangle$ and let $b_n = a_n \vartriangle a$,
for each $n \in \omega$.

The sequence $\langle b_n : n \in \omega \rangle$ converges
to $\bold 0$. We claim that $\langle b_n : 
n \in \omega \rangle$ is a name for a finite set. To see this,
let $\langle \tilde b_n : n \in \omega \rangle$ be a decreasing
sequence with $b_n \leq \tilde b_n$ and $\bigwedge_n \tilde b_n 
= \bold 0$. If $G$ is a generic filter on $B$, then only finitely
many $\tilde b_n$'s can be in $G$, hence $\langle \tilde b_n : 
n \in \omega \rangle$ is a name for a finite set, and so is 
$\langle b_n : n \in \omega \rangle$.

Now the constant sequence $\langle a : n \in \omega \rangle$
is a name for either $\omega$ (with Boolean value $a$) or 
$\emptyset$ (with value $-a$), and since $a_n = a \vartriangle b_n$,
the $B$-valued real $\langle a_n : n \in \omega \rangle$ is the 
symmetric difference of either $\omega$ or $\emptyset$ and a finite
set. Hence the real $\langle a_n : n \in \omega \rangle$ is either
finite or cofinite, and hence not independent.

Conversely, let $\langle a_n : n \in \omega \rangle$ be a sequence
that has no convergent subsequence. We shall produce a name for
an independent real (or rather independent with nonzero Boolean
value).

First we claim that $\langle a_n : n \in \omega \rangle$ has a 
subsequence $\langle c_n : n \in \omega \rangle$ with the 
property that $\limsup x_n = \limsup c_n$ for every subsequence
$\langle x_n : n \in \omega \rangle$ of $\langle c_n : 
n \in \omega \rangle$. This is proved as follows: let 
$a_n^0 = a_n$. Suppose that $\langle a_n^0 : n \in \omega \rangle$
has a subsequence $\langle a_n^1 : n \in \omega \rangle$
with $\limsup a_n^0 > \limsup a_n^1$, which has a subsequence
$\langle a_n^2 : n \in \omega \rangle$ with $\limsup a_n^1
 > \limsup a_n^2$ and so on. At limit stages we produce a 
subsequence by diagonalization. Since $B$ satisfies $ccc$,
the process stops after countably many steps, and we obtain
$\langle c_n : n \in \omega \rangle$ with the desired 
property.

We repeat this argument for $\liminf$, and so we may assume
that $\liminf x_n = \liminf c_n$ for every subsequence
$\langle x_n : n \in \omega \rangle$. Since $\langle c_n :
 n \in \omega \rangle$ is not convergent we have $\liminf c_n
< \limsup c_n$. Let $c = \liminf c_n$ and let $b_n = c_n - c$,
for each $n \in \omega$. We have $\liminf x_n = \bold 0$ and
$\limsup x_n = u > \bold 0$ for every subsequence $\langle x_n 
: n \in \omega \rangle$ of $\langle b_n : n \in \omega \rangle$.

It follows that for every $d$ with $\bold 0 < d \leq u$, 
$d \wedge b_n \neq \bold 0$ and $d - b_n \neq \bold 0$ for
all but finitely many $n$: otherwise, there is a subsequence
$\langle x_n : n \in \omega \rangle$ such that either $x_n \geq d$
for all $n$ or $x_n \leq -d$ for all $n$, contradicting
$\liminf x_n = \bold 0$ and $\limsup x_n = u$. We claim that
$u$ forces that $\langle b_n : n \in \omega \rangle$ is an
independent real. Let $G$ be a generic on $B$ with $u \in G$, 
and let $X = \{n:b_n \in G \}$. If $X$ is not independent then
there exists an infinite $A \subset \omega$ (in $V$) such
that either $A \subset X$ or $A \cap X = \emptyset$. In the
former case let $d = u \wedge \bigwedge_{n \in A} b_n$; in
the latter, let $d = \bigwedge_{n \in A} {u -b_n}$. In either case, 
$d \in G$ and so $\bold 0 < d \leq u$ and either $d \leq b_n$
or $d \leq -b_n$ for infinitely many $n$, a contradiction.
\hfill $\square$
\enddemo

\remark{Remarks} \roster \runinitem "(i)" In \cite{BFH}
there is another, algebraic, characterization of Boolean algebras
$B$ that do not add independent reals: the existence of an 
{\it almost regular} embedding of the Cantor algebra (i.e. the 
countable atomless Boolean algebra) into~$B$.\endroster
\roster
\runinitem "(ii)" In sequential $T_1$ spaces, sequential
compactness is equivalent to countable compactness
(every countable open cover has a finite subcover).
Therefore a complete $ccc$ Boolean algebra $B$ does
not adds independent reals if and only if $(B,\tau_s)$
is countably compact.\endroster
\roster
\runinitem "(iii)" If $B = \Cal P (\omega)$ then $(B,\tau_s)$
is a compact space. We don't know whether there exists
an atomless Boolean algebra $B$ such that $(B,\tau_s)$
is compact. Note that $B$ has to be $ccc$.
\endroster
\endremark

%.............................5. examples....................

\subheading{5. Examples} We present three examples of complete
$ccc$ Boolean algebras that do not add independent reals. All
three examples are only consistent, not in ZFC. 
\medskip

\noindent {\bf 5.1} The first example is a complete Boolean algebra 
$B = \Cal P (\kappa) / I$
where $I$ is certain $\sigma$-saturated ideal on $\kappa$.
The example is due
to G\l \'owczy{\'n}ski who showed in \cite{G} that $B$ is weakly 
distributive, countably generated and does not carry 
a strictly positive Maharam submeasure. We show that $B$ does 
not add independent reals.

We use the known properties of the sequential topology on
the power set algebra $\Cal P (\kappa)$, cf. \cite{BGJ}:
\roster
\item "(i)" $(\Cal P (\kappa),\tau_s)$ is Hausdorff
\item "(ii)" $(\Cal P (\kappa),\tau_s)$ is regular if and only if
$\kappa = \omega$
\item "(iii)" $(\Cal P (\kappa),\tau_s)$ is Fr\' echet if and
only if $\kappa < \goth b$
\item "(iv)" $(\Cal P (\kappa),\tau_s)$ is sequentially compact
if and only if $\kappa < \goth s$ ($\goth s$ is the splitting number)
\endroster

If $I$ is $\sigma$-ideal on $\Cal P (\kappa)$ then it is a closed
subset in $(\Cal P (\kappa),\tau_s)$ and the sequential topology
on $B = \Cal P (\kappa) / I$ is the quotient topology of $\tau_s$.
Hence if $(\Cal P (\kappa),\tau_s)$ is Fr\' echet (or sequentially
compact) then so is $(B, \tau_s)$. 

Now let $\kappa$ be a measurable cardinal and let $V[G]$ be a
generic extension of $V$ by $ccc$ forcing
in which Martin's Axiom holds and $2^{\aleph_0} > \kappa$. The
measure on $\kappa$ in $V$ generates a nonprincipal 
$\sigma$-saturated $\sigma$-ideal $I$ on $\kappa$ in $V[G]$.
Let $B = \Cal P (\kappa) /I$. $B$ is an atomless complete
$ccc$ Boolean algebra, and since $\goth b = \goth s =
2^{\aleph_0} > \kappa$ (from MA) it follows that $(B,\tau_s)$ is
Fr\' echet and sequentially compact. Hence $B$ is weakly
distributive and does not add independent reals.

\medskip

\noindent {\bf 5.2.} The second example is Jensen's forcing 
\cite{J} that
produces a minimal nonconstructible real. The corresponding
Boolean algebra (constructed in $L$) is $ccc$ and weakly
distributive. A slight modification of Jensen's construction 
guarantees that the forcing does not add independent reals.

We proceed under the assumption of $V=L$, and assume that
$\{S_\alpha : \alpha < \omega_1 \}$ is a diamond sequence,
namely such that for every $A \subset L_{\omega_1}$, the
set $\{\alpha < \omega_1 : A \cap L_\alpha = S_\alpha \}$
is stationary. 

Let $Seq=\{0,1\}^{<\omega}$ be the set of all finite $0\!-\!1$ sequences.
We shall construct a forcing notion $P$
consisting of perfect trees $T \subset Seq$; the ordering
of $P$ is by inclusion. $P$ will be the union of a continuous
$\omega_1$-sequence of countable sets $$P_0 \subset P_1
\subset \dots \subset P_\alpha \subset \dots \quad
\alpha < \omega_1$$ where every $P_\alpha$ is closed under
taking restrictions $T \upharpoonright s = \{t \in T :
t \subset s \text{ or } t \subset s \}$ (where $s \in T$).
Let $P_0 = \{ T_0 \upharpoonright s : s \in Seq \}$ where
$T_0=Seq$. At limit stages, $P_\alpha = \bigcup_{\beta <
\alpha} P_\beta$.

We now describe the construction of $P_{\alpha + 1}$
from $P_\alpha$. Let ${\Cal X}_\alpha$ be the set of all
$X$ such that for some $\beta \leq \alpha$, $X = S_\beta$
and is a predense set in $P_\beta$, along with all
$Q_\beta = P_{\beta + 1} - P_\beta$ for $\beta < \alpha$.
The inductive condition is that each $X \in {\Cal X}_\alpha$
is predense in $P_\alpha$. (This inductive condition
remains true at the limit stages.) Enumerate the countable
set ${\Cal X}_{\alpha}$ so that each $X$ occurs infinitely
often in the enumeration: ${\Cal X}_\alpha = \{X_n^\alpha :
n \in \omega \}$. For each $p \in P_\alpha$ we construct
a perfect tree $T = T(\alpha,p) \subset p$ and then let
$Q_\alpha = \{T(\alpha,p) \upharpoonright s : p \in P_\alpha,
\ s \in T(\alpha,p) \}$ and $P_{\alpha + 1} = P_\alpha \cup 
Q_\alpha$. The tree $T$ will be the fusion of a collection
$\{p_\sigma : \sigma \in Seq \}$ where each $p_\sigma$ is
in $P_\alpha$ and:
\roster
\item "(i)" $p_\emptyset \subset p$,
\item "(ii)" $p_{\sigma^{\smallfrown} 0}$ and 
$p_{\sigma^{\smallfrown} 1}$ are both stronger then $p_\sigma$
and have incompatible stems,
\item "(iii)" $T = \bigcap_{n \in \omega} \bigcup_{|\sigma| = n} 
p_\sigma.$
\endroster

If it is not the case that $S_\alpha$ is a $P_\alpha$-name 
for a subset of $\omega$, we let $p_\emptyset = p$,
and for each $\sigma \in Seq$, if $|\sigma|=n$, find 
$p_{\sigma^{\smallfrown} 0}$ and $p_{\sigma^{\smallfrown} 1}$
in $P_\alpha$ that satisfy $(ii)$ such that
\roster
\item "(iv)" both $p_{\sigma^{\smallfrown} 0}$ and 
$p_{\sigma^{\smallfrown} 1}$ are stronger than some
$x \in X_n^\alpha$.
\endroster

 If $S_\alpha = \dot A$ is a $P_\alpha$ name for a subset
of $\omega$ and some $q \subset p$ forces infinitely many $n$ into
$\dot A$, we let $p_\emptyset = q$ and again find $p_\sigma$,
$\sigma \in Seq$, that satisfy $(ii)$ and $(iv)$.

Thus assume 

$$\forall q \subset p \ \exists k \ \forall n \geq k \ \exists r 
\subset q \
r \Vdash n \notin \dot A$$
[We wish to point out that $q$ and $r$ ranges over $P_\alpha$,
and the forcing relation $\Vdash$ refers to the forcing $P_\alpha$.]

At stage $n$, first find for each $\sigma$ with $|\sigma| = n$ 
conditions ${\bar p}_{\sigma^{\smallfrown} 0}$ and 
${\bar p}_{\sigma^{\smallfrown} 1}$ that satisfy $(ii)$ and
$(iv)$, and then find some $a_n \in \omega$ sufficiently large
so that $a_n > a_{n-1}$ and for each $\sigma$ with $|\sigma| =n$
there exist conditions $p_{\sigma^{\smallfrown} 0} \subset 
{\bar p}_{\sigma^{\smallfrown} 0}$ and $p_{\sigma^{\smallfrown} 1} 
\subset {\bar p}_{\sigma^{\smallfrown} 1}$ that all force
$a_n \notin \dot A$. Then let $Y_p^\alpha = \{a_n : n \in \omega\}$,
let $T$ be the fusion of $\{p_\sigma : \sigma \in seq\}$, and finally,
$Q_\alpha = \{T(\alpha,p)\upharpoonright s : p \in P_\alpha, \ 
s \in T(\alpha,p) \}$.

We claim that each $X \in {\Cal X}_\alpha$ is predense in 
$P_{\alpha+1} = P_\alpha \cup S_\alpha$: Let $q\in P_{\alpha+1}$.
As $X$ is predense in $P_\alpha$, we may assume that 
$q \in Q_\alpha$, $q = T(\alpha,p) \upharpoonright s$.

Let $\{p_\sigma : \sigma \in Seq \}$ be the fusion collection
for $T(\alpha,p)$, and let $n > |s|$ be such that
$X = X_n^\alpha$. There exists some $\sigma$ with $|\sigma| =
n+1$ such that $t \supset s$ where $t$ is the stem of
$p_\sigma$, and by $(iv)$, $p_\sigma$ is stronger then some
$x \in X$. Hence $q \upharpoonright t < p_\sigma < x$, and
so $q$ and  $x$ are compatible in $P_{\alpha+1}$.

It follows that every $X \in {\Cal X}_\alpha$ is predense in
every $P_\beta,\ \beta \geq \alpha$, and therefore in $P$.

Now it follows that $P$ satisfies $ccc$: Let $X$ be a maximal
antichain in $P$. For a closed unbounded set of $\alpha$'s,
$X \cap \alpha$ is a maximal antichain in $P_\alpha$. Therefore
there exists an $\alpha$ such that $X \cap \alpha = S_\alpha$
is a maximal antichain in $P_\alpha$ and hence $X \cap \alpha
\in {\Cal X}_\alpha$. Thus $X \cap \alpha$ is predense, and 
hence a maximal antichain, in $P$, and so $X = X \cap \alpha$
and $X$ is countable.

It is well known that forcing with perfect trees does not add
unbounded reals and so $P$ is weakly distributive. We shall now prove
that $P$ does not add independent reals. Let $\dot A$ be a name for
a subset of $\omega$, and let $p$ be a condition. We prove that
there exists a stronger condition $q$ and some infinite set 
$Y \subset \omega$ such that either $q \Vdash Y \subset \dot A$
or $q \Vdash Y \cap \dot A = \emptyset$. Thus assume that there
is no $q \subset p$ that forces infinitely many $n$ into $\dot A$.

For each $n$, let $X_n$ be a maximal antichain whose members all
decide $n \in \dot A$, and let $\gamma$ be large enough so that
$X_n \subset P_\gamma$ for each $n$. Note that for every 
$\alpha \geq \gamma$ and every $q \in P_\alpha$, 
$q \Vdash n \in \dot A$ has the same meaning in $P_\gamma$ as in $P$.
So let $\alpha \geq \gamma$ be such that $\dot A = S_\alpha$.
Let $\{p_\sigma : \sigma \in Seq \}$ be the fusion collection for 
$T(\alpha,p)$. As there is no $q \subset p$ in $P_\alpha$ that
forces infinitely many $n$ into $\dot A$, there exists an infinite
set $Y = Y_p^\alpha = \{a_n:n \in \omega \}$ such that for every
$n$ and every $\sigma$ with $|\sigma| = n+1$, $p_\sigma \Vdash
a_n \notin \dot A$. It follows that the condition $T(\alpha,p)$
forces $Y \cap \dot A = \emptyset$. Hence $P$ does not add
independent reals.

\medskip

\noindent {\bf 5.3.} The third example is a complete $ccc$ 
Boolean algebra that is 
not weakly distributive and does not add independent reals.
Again, we work under the assumption of $V=L$. While the previous
example is a $ccc$ version of Sacks forcing \cite{S}, this 
example is a $ccc$ variant of Miller's forcing \cite{Mi} with
superperfect trees. We show that Miller's argument for the
absence of independent reals can be used in the context of this
Jensen~-~style construction.

Let $\{S_\alpha : \alpha \in \omega_1 \}$ be a diamond sequence
for $L_{\omega_1}$, and let $Seq = \omega^{<\omega}$. Forcing
conditions will be superperfect trees $T \subset Seq$ and $P$ will
be constructed via a continuous sequence
$$P_0 \subset P_1 \subset \dots \subset P_\alpha \subset \dots, \quad
\alpha \in \omega_1$$
of countable sets closed under restrictions, with $P_0 = \{ Seq
\upharpoonright s : s \in Seq \}$ and $P_\alpha = \bigcup_{\beta <
\alpha} P_\beta$ for limit $\alpha$'s.

At stage $\alpha$ of the construction, let ${\Cal X}_\alpha$ be
the set of all $S_\beta,\ \beta \leq \alpha$, that are predense
in $P_\beta$ and all $Q_\beta = P_{\beta +1} - P_\beta$, 
$\beta <\alpha$. By induction, each $X \in {\Cal X}_\alpha$ is
predense in $P_\alpha$. Enumerate ${\Cal X}_\alpha$ so that
each $X$ occurs infinitely often, ${\Cal X}_\alpha = 
\{X_n^\alpha : n \in \omega\}$. For each $p \in P_\alpha$ 
we construct a superperfect tree $T(\alpha,p) \subset p$ and
let $Q_\alpha =\{T(\alpha,p)\upharpoonright s : p \in P_\alpha, \
s \in T(\alpha,p)\}$, and $P_{\alpha+1} = P_\alpha \cup Q_\alpha$.

Assume that $S_\alpha = \dot A$ is a name for a subset of $\omega$;
along with $T(\alpha,p)$ we construct an infinite set 
$Y_p^\alpha = \{a_n : n \in \omega \}$ such that (under the right
circumstances) the condition $T(\alpha,p)$ will force (in $P$)
either $Y_p^\alpha \subset \dot A$ or $Y_p^\alpha \cap \dot A = 
\emptyset$.
At the same time, $T(\alpha,p)$ will be compatible with every
$X \in {\Cal X}_\alpha$, to guarantee that $X$ remains predense
in $P_{\alpha+1}$. [If $S_\alpha$ is not a name for a subset of
$\omega$ then we only handle second requirement at stage $\alpha$.]

We recall that $s \in T$ is a {\it splitting node} if 
$s^{\smallfrown}k \in T$ for infinitely many $k$, and that $s$ is an
$n^{\text{th}}$ splitting node if moreover $|\{t: t\subset s \ 
\text{is  a  splitting node} \}| = n$.

\smallskip

\noindent {\it Step 1.} Let $U$ be a nonprincipal ultrafilter on
$\omega$. We construct a superperfect tree $T \subset p$, and for
each splitting node $s \in T$ a set $A_s \subset \omega$, and for
each successor $s^{\smallfrown}k$ of $s$ in $T$ a condition 
$p_{s^{\smallfrown}k} \in P_\alpha$ such that $s^{\smallfrown}k
\subset \text{stem}(p_{s^{\smallfrown}k})$ and 
\roster
\item "(1)" $T \upharpoonright (s^{\smallfrown}k) \subset p_s$,
\item "(2)" if $s$ is an $n^{\text{th}}$ splitting node then
$p_{s^{\smallfrown}k} \subset x$ for some $x \in X_n^\alpha$,
\item "(3)" either $\forall s \ A_s \in U$ or $\forall s \ 
- A_s  \in U$,
\item "(4)" $\forall m \ \exists k_0 \ \forall k \geq k_0$ if 
$s^{\smallfrown}k
\in T$ then $p_{s^{\smallfrown}k} \Vdash \dot A \cap m = A_s \cap m$.
\endroster

To construct $T$, let $T_0 = p$ and let $s = \text{stem}(p)$ be the
1${}^{\text{st}}$ splitting node of $T_0$. For each successor 
$s^{\smallfrown}k$ of $s$ in $T_0$, find $p_{s^{\smallfrown}k}$ with
$s^{\smallfrown}k \subset \text{stem}(p_{s^{\smallfrown}k})$ such
that $(2)$ holds for $X_0^\alpha$ and such that $p_{s^{\smallfrown}k}$
decides $\dot A \cap k$. Then thin out the successor $s^{\smallfrown}k$
successively so that $0 \in \dot A$ is decided the same way by all,
$1 \in \dot A$ is decided the same way by all starting with the 
second one, $2 \in \dot A$ by all starting with the third one etc.
When finished, let $A_s = \{ m : \ \text{eventually all} \
p_{s^{\smallfrown}k} \ \text{force}\  m \in \dot A\}$, and let 
$T_1 = \bigcup \{p_{s^{\smallfrown}k} : s^{\smallfrown}k \ \text{are 
the retained successors of}\ s\}$.

Next consider all 2${}^{\text{nd}}$ splitting nodes $s$ of $T_1$ and
repeat the construction of $T_1$ from $T_0$, using $X_1^\alpha$.
Repeating this $\omega$ times, we get trees $T_n,\ n \in \omega$, 
and let $\tilde T = \bigcap_{n=0}^\infty T_n$; $\tilde T$ is a
superperfect tree. Let $\tilde T = C_1 \cup C_2$ where
$C_1 = \{s: A_s \in U \}$ and $C_2 = \{s: - A_s \in U \}$.
Either $C_1$ or $C_2$ contains a superperfect tree $T$. The tree
$T$ satisfies $(1)$~-~$(4)$.

\smallskip

\noindent {\it Step 2.} Assume that $\forall s \ A_s \in U$
(the argument is similar in the opposite case). Let $T_0 = T$,
$F_0 =\{s_0 \}$ where $s_0 = \text{stem}(T_0)$ and let 
$a_0 \in A_s$. All but finitely many successors $s^{\smallfrown}k$
in $T_0$ have the property that $p_{s^{\smallfrown}k} \Vdash
a_0 \in \dot A$, and so remove the finitely many successors,
resulting in a tree $T_1 \subset T_0$ with stem $s_0$.

Let $F_1 \subset T_1$ be the $\{s_0,s_1\}$ where $s_1$ is the 
leftmost 2${}^{\text{nd}}$ splitting node of $T_1$, and let $a_1 > a_0$
be such that $a_1 \in A_{s_0} \cap A_{s_1}$. (The set is nonempty 
because in $U$.) Remove finitely many successors of $s_1$ so that
for all the remaining ones, $p_{s_1^{\smallfrown}k} \Vdash a_1 \in \dot 
A$.
Also remove finitely many successors of $s_0$, with the same result.
(There is no problem with the successors of $s_0$ below $s_1$, 
because every node above it is either below or above $s_1$.)
The resulting $T_2$ is such that $F_1 \subset T_2 \subset T_1$.

Now let $F_2 \subset T_2$ be the finite set $F_2 \supset F_1$ that
is obtained by adding to $F_1$ the leftmost 3${}^{\text{rd}}$
splitting node $s_2$ above $s_1$ and the second leftmost
2${}^{\text{nd}}$ splitting node $s_3$ (above $s_0$). Let 
$a_2 > a_1$ be such that $a_2 \in \bigcup_{s \in F} A_s$.
Remove finitely many successors of $s_3$, finitely many successors 
of $s_2$, finitely many successors of $s_1$ and finitely many 
successors 
of $s_0$ (in that order)
so that all the remaining $p_{s_i^{\smallfrown}k}$ force 
$a_2 \in \dot A$. This produces $T_3$, and $F_2 \subset T_3 
\subset T_2$. We continue in this fashion, and let 
$T(\alpha,p) = \bigcap_{n=0}^\infty T_n$, $Y_p^\alpha =
\{a_n : n \in \omega \}$. This completes the construction of 
$P_{\alpha +1}$.

As in the second example, the forcing $P$ satisfies $ccc$ because
every $X \in {\Cal X}_\alpha$ is predense in every $P_\beta$, 
$\beta \geq \alpha$. It is well known that forcing with 
superperfect trees does add an unbounded real and so $P$
is not weakly distributive.

If $\dot A$ is a name for a subset of $\omega$ and $p \in P$,
then for some sufficiently large $\alpha$, $\dot A = S_\alpha$
and $q \Vdash n \in \dot A$ has the same meaning in $P_\alpha$
as in $P$. The construction of $T(\alpha,p)$ and $Y_p^\alpha$
yields that the condition $T(\alpha,p)$ either forces
$Y_p^\alpha \subset \dot A$ or forces $Y_p^\alpha \cap \dot A =
\emptyset$.

%.............................references.............................

\enddocument